\documentclass[12pt]{article}
\usepackage{graphics, latexsym}
\usepackage{graphicx}
\usepackage{setspace}
\usepackage{amssymb}
\usepackage{amsmath}
\usepackage{float}
\usepackage{authblk}
\usepackage{booktabs}
\usepackage{threeparttable}
\usepackage{multirow}
\usepackage{rotating}
\begin{document}
\title{ Homogeneity Test of Several High-Dimensional Covariance Matrices for Stationary Processes under Non-normality }
\author{ABDULLAH QAYED$ ^{1} $ and DONG HAN$ ^{1} $\\
Department of Statistics$ ^{1} $, Shanghai Jiao Tong University, China.}

\maketitle
\begin{abstract}
 This article presents a homogeneity test for testing the equality of several high-dimensional covariance matrices for stationary processes with ignoring the assumption of normality. We give the asymptotic distribution of the proposed test. The simulation illustrates that the proposed test has perfect performance. Moreover, the power of the test can approach any high probability uniformly on a set of covariance matrices.\end{abstract}

\renewcommand{\thefootnote}{\fnsymbol{footnote}}
\footnotetext{ $ ^{1} $ Supported by National Basic Research Program of China
(973 Program, 2015CB856004) and the National Natural Science Foundation of China (11531001). $\,\, MSC\, 2010\,subject\,classification $. Primary 62L10; Secondary 62L15}
\footnotetext{CONTACT: ABDULLAH QAYED, email: qayed@sjtu.edu.cn - Department of Statistics-Shanghai Jiao Tong University. Shanghai, Minhang District. Shanghai 200240,China.}

\section{Introduction}

Many statistical methods in multivariate analysis require the assumption of homogeneous covariance matrices; for instance,  multivariate analysis of variance, Hotelling $T^2$-test, discriminant function analysis and Mahalanobis' generalized distance studies,  in addition to the statistical models in repeated measures, longitudinal studies and multilevel regression analysis and so on. Moreover, the inference results of many applied sciences are depend on  testing the covariance matrices, as in biometrics in the study of homogeneity of covariance matrices in infraspecific paleontologic and neontologic material and, if possible, by comparing with material known to be homochronous.
On other hand, with the expansion of exact procedures for measurement theory in many applied sciences such as biology, life and other sciences, data of many different research phenomena may not be distributed normally, as in the exponential distribution for growth, microorganisms, pandemics, cancer cells and compound interest in economic surveys,... etc. Other example, in genomics, Torrente et al (2019) investigated the prevalence of genes with expression distributions, they found that there are less than $50\%$ of all genes were Normally-distributed, with other distributions including Gamma, Bimodal, Cauchy, and Log-normal were represented.  In such a case, when studying a gene's distribution shape to test heterogeneity of data, the variety of non-normal data distributions is often overlooked when using homogeneity tests that assume all data groups follow the normal distribution.

In this paper, we consider homogeneity tests of covariance matrices with high dimensional data when $p$ can be much larger than the sample size $n_i$.
Recently, many homogeneity tests have been developed by many authors to address the deficiencies of classical tests caused by high-dimensional $p$ and corresponding small sample sizes, and these tests are based on the trace of sample covariance matrices, as in sphericity hypothesis tests proposed by Jiang \&  Chen et al (2010), Qiu and Chen (2012), Ledoit and Wolf (2002), Srivastava (2005), Cai et al (2013), Cai and Ma (2013),  Li and Chen (2014) ,Peng et al (2016), Ishii et al (2019) and Zhang et al (2020). Despite the classical homogeneity tests (as in Box (1949), Jiang et al (2012) and Jiang \& Yang (2013)) are sensitive to departure the normal distribution, almost recently developed tests also assume the normal distribution of data. This raises questions about the accuracy of the various homogeneity tests with relatively small sample sizes of groups or when assumptions of normality and large samples are not met.

Moreover, for testing the equality of several covariance matrices under the hypotheses
\begin{align}\label{(1)}
 H_{0}:\Sigma_{1}=...=\Sigma_{k}=\Sigma  \,\,\, vs\,\,\,  H_{1}: \Sigma_{i} \neq \Sigma_{j}, \text{ for at least one pair $ (i,j),\,\,\, i\neq j $},
\end{align}
where $\Sigma_{i}$ denotes the the $p\times p$ covariance matrices of the $i$th population with $p$-dimensional multivariate distribution for $1\leq i\leq k$, where $k\geq 2$. Schott (2007) presented the $J_k$ test based on the trace of sample covariance matrices and proved that the asymptotic distribution of $J_k$ is normal $N(0, 1)$ under the null hypothesis $H_0$ in (1) with the setting of $p/n_i \to b_i \in (0, \, \infty)$ as $n_i, p \to \infty$ for $1\leq i\leq k$, where $n_i$ is the number of samples for $i$th population.   Srivastava and Yanagihara (2010)  proposed two tests $T^{2}_k$ and  $Q^{2}_k$ based on the trace of sample covariance matrices for testing $  H_{0} $ and showed that the asymptotic null distributions of both testes are  $ \chi^{2}_{k-1} $ in the setting of $n_i=O(p^{\delta}), \delta >1/2$ as $n_i, p \to \infty$ for $1\leq i\leq k$.
However, all tests $J_k$, $T^{2}_k$ and  $Q^{2}_k$ assume the normal distribution of data.
Based on the $U$-statistics,  Ahmad (2016) and Zhong et al (2019) presented tests $T_g$ and $\hat{D}_{nt}$ respectively for high-dimensional data under non-normality and  high-dimensional longitudinal data,  where both tests have asymptotic null normal distribution and do not need a specific relationship between $p$ and $n_i$ as $p=p(n_i), n_i \to \infty$. However, all five tests $J_k$, $T^{2}_k$, $Q^{2}_k$, $T_g$ and $\hat{D}_{nt}$ above need to assume that $0<a_i=\lim_{p\to \infty}tr(\Sigma^j)/p<\infty, \, j=1, 2, 3, 4$, or $0<tr(\Sigma^2)/p^2 = \delta_1 <\infty$ for all $p$ and $inf_{p}\{tr(\Sigma^4)/[tr(\Sigma^2)]^2\}>0$.

The research in this paper is motivated by the conventional likelihood ratio $M$ test in Box (1949), which established the asymptotic $\chi^2$ distribution of the test statistic for fixed-dimensional normally distributed random vectors, when  $n_i > p$,  $1\leq i\leq k$. In this paper, a modified $M$ test, $L_k$, is proposed to test the null hypothesis in (1) without the normality assumption. We don't assume the conditions about the relationship between the traces of powers of covariance matrix and $p$ listed above for testing the null hypothesis $H_0$.
However,  we assume that the relationship between $p$ and $n_i$ satisfies $p\geq c(n_{max})^{3r/(r-2)} $ for some positive constant $c$, where $n_{max}=\max\{n_1, n_2, ...,n_k\}$ and $r>2$. Hence, the test $L_k$ accommodates the large $p$, small $n_{max}$ situations.

This paper is organized as follows. Section 2 introduces the basic data structure of the G-stationary process and its associated properties, and a modified $M$ test. In section 3, we present the main results of this paper, we show that the asymptotic null distribution of $L_k$, as $n_{min} \to \infty$, is $ \chi^{2}_{k-1}$ under the null hypothesis $H_0$ in (1) followed by the procedures of $L_k$ test. Reports of numerical simulation are given in section 4. An empirical study using time-series data to  extract clinical information from speech signals is presented in section 5. The paper concludes in section 6. The proofs of Theorems 3.1 and 3.2, and the examples results are deferred to the Appendix.

\section{The G-stationary process and a modified $M$ test}

We first introduce the G-stationary process and its associated properties, then present a modified $M$ test.

A dependent stationary process is called as the G-stationary process if it has the following representation
\begin{align*}
X_j=G(...,\varepsilon_{j-1}, \varepsilon_{j}, \varepsilon_{j+1},...), \,\,\,\,\,\,\,          j\in Z
\end{align*}
and satisfies the conditions of Theorem 2.1 in Berkes, Liu and Wu (2014), where $\varepsilon_{i}, i\in Z$, are i.i.d. random variables,  $G: R^Z\rightarrow R$ is a measurable function, $Z$ and $R$ denote the set of integers and the set of real numbers, respectively. In fact, under some appropriate conditions, the G-stationary process can include a large class of linear and nonlinear processes, such as the functionals of linear processes, bilinear models, GARCH processes, generalized random coefficient autoregressive models, nonlinear AR models etc., see  Liu and Lin (2009), Berkes, Liu and Wu (2014).

Assume that $X_j$ has mean zero, $\textbf{E}(|X_j|^r)<\infty$ for $r>2$, with covariance function $\gamma(j)=\textbf{E}(X_0X_j), j\in Z$. It is clear that the covariance matrix $\Sigma_p=(\sigma_{ij})_{p\times p}$ of $(X_1, X_2,...,X_p)$ satisfies $\Sigma_p=(\sigma_{ij})_{p\times p}=\Sigma_p=(\gamma(j-i))_{p\times p}$.  Let $\mathbb{S}_p=\sum_{j=1}^pX_j$. The most important result for the G-stationary process is that it satisfies the strong invariance principle, i.e., it can arrive the optimal KMT approximation (see Berkes, Liu and Wu 2014), that is,  there exists a richer probability space such that
\begin{align}\label{(2)}
\mathbb{S}_p-\sigma\mathbb{B}(p)=o(p^{1/r}), \,\,\,\, a.s.
\end{align}
where $\mathbb{B}(t)$ is a standard Brownian motion  and $\sigma^2=\sum_{j\in Z}\gamma(j)$. From (2) it follows that
\begin{align}\label{(3)}
\frac{\mathbb{S}_p}{\sigma\sqrt{p}} \Rightarrow N(0, 1) \,\,\,\, \text{ and } \,\,\,\,\, \sigma^2=\lim_{p\to \infty}\frac{\textbf{E}(\mathbb{S}^2_p)}{p}.
\end{align}

In this paper, our main purpose is to present a test for testing the equality of $k$ covariance matrices $\Sigma_1, \Sigma_{2},...,\Sigma_{k}$ which correspond to $k$ G-stationary processes, $X^{(i)}_j, 1\leq j \leq p$ for $1\leq i\leq k$, respectively, where $\Sigma_{i}=(\sigma^{(i)}_{jj'})_{p\times p}=(\gamma^{(i)}(j-j'))_{p\times p}$ and $\gamma^{(i)}(j)=\textbf{E}(X^{(i)}_0X^{(i)}_j), j\in Z$ for $1\leq i\leq k$. The $k$ G-stationary processes can be considered as $k$ $p$-dimensional random vectors, $X^{(i)}=(X^{(i)}_1, X^{(i)}_2, ..., X^{(i)}_p), \, 1\leq i\leq k$. For each $i$, we take $n_i$ i.i.d. random samples (random vectors)
\begin{align*}
X^{(i)}_l=(X^{(i)}_{1l}, X^{(i)}_{2l}, ..., X^{(i)}_{pl}), \,\,\,\,\,\,\,\, 1\leq l\leq n_i
\end{align*}
to estimate $\Sigma_{i}$ for $1\leq i\leq k$. Sometimes, $p$ is much bigger than $n_i$. For example, consider whether the fluctuation of hypertension in half a year ( 180 days ) in patients with hypertension in $k$ countries or regions is consistent. We may randomly select 100 patients with hypertension in every country and wear a smart health bracelet to every patient with hypertension, which can monitor the blood pressure profile of the patient every minute. Here, $n_1=n_2=...=n_k=100$ and $p=60\times 24\times 180=259200$.

Let $ \textit{\textbf{S}}_{i} $ be the sample covariance matrix for $\Sigma_{i}, 1\leq i\leq k.$  Note that $ \textit{\textbf{S}}_{i} $ can be written as
\begin{align}\label{4}
\textit{\textbf{S}}_{i}=  \frac{1}{n_{i}-1} \sum_{l=1}^{n_{i}} \left(\begin{array}{c}
X^{(i)}_{1l} - \overline{X}^{(i)}_{1}\\
X^{(i)}_{2l} - \overline{X}^{(i)}_{2}\\
...\\
X^{(i)}_{pl} - \overline{X}^{(i)}_{p}
\end{array}\right)\,\,\, \left(\begin{array}{cccc}
X^{(i)}_{1l} - \overline{X}^{(i)}_{1},\, X^{(i)}_{2l} - \overline{X}^{(i)}_{2},\, ...,\, X^{(i)}_{pl} - \overline{X}^{(i)}_{p}
\end{array}\right)
\end{align}
and therefore
\begin{align}\label{5}
 \textit{\textbf{1}}\textit{\textbf{S}}_{i}\textit{\textbf{1}}^{T} = \frac{1}{n_i-1} \sum_{l=1}^{n_{i}}\Big(\sum_{j=1}^{p}(X^{(i)}_{jl}-\overline{X}^{(i)}_{j})\Big)^2\geq 0,
\end{align}
where $ \overline{X}^{(i)}_{j}=\sum_{l=1}^{n_{i}} X^{(i)}_{jl}/n_i $ for $1\leq j\leq p$ and $ \textit{\textbf{1}}=(1, 1,...,1)$ is $p$-dimensional vector with $1$ as the component.
How to guarantee  $\textit{\textbf{1}}\textit{\textbf{S}}_{i}\textit{\textbf{1}}^{T} >0, \, a.s.$ ? In fact,  if  $\mathbb{S}^{(i)}_{pl}:=\sum_{j=1}^pX^{(i)}_{jl}$ is a continuous random variable for each $i, 1\leq i\leq k$, then
\begin{align*}
\mathbb{S}^{(i)}_{pl}-\frac{1}{n_i}\sum_{l=1}^{n_i}\mathbb{S}^{(i)}_{pl}=(1-\frac{1}{n_i})\mathbb{S}^{(i)}_{pl}-\frac{1}{n_i}\sum_{l'\neq l}\mathbb{S}^{(i)}_{pl'} \neq 0,\,\, a.s.
\end{align*}
for $n_i>1, 1\leq i\leq k$,  since $\mathbb{S}^{(i)}_{pl},  1\leq l\leq n_i$,  are i.i.d. continuous random variables for each $i, 1\leq i\leq k$. This means that
\begin{align*}
(n_i-1)\textit{\textbf{1}}\textit{\textbf{S}}_{i}\textit{\textbf{1}}^{T}=\sum_{l=1}^{n_{i}}\Big(\mathbb{S}^{(i)}_{pl}-\frac{1}{n_i}\sum_{l=1}^{n_i}\mathbb{S}^{(i)}_{pl})\Big)^2 \, >0, \,\, a.s.
\end{align*}
for $n_i>1, 1\leq i\leq k$. For general case, we can prove that $\textit{\textbf{1}}\textit{\textbf{S}}_{i}\textit{\textbf{1}}^{T}/p \,\, >0,  \,\, a.s.$ as $p \to \infty $ in Theorem 3.1 in the next section.

The often used test for testing the equality of several covariance matrices is the Box's $ M $ (1949) test which can be written in the following
\begin{align}\label{6}
M=& (n-k)log |\textit{\textbf{S}}|-\sum_{i=1}^{k} (n_{i}-1)log |\textit{\textbf{S}}_{i}|,
\end{align}
and $\textit{\textbf{S}}=\sum_{i=1}^{k} (n_{i}-1)\textit{\textbf{S}}_{i}/(n-k)$ and $n-k=\sum_{i=1}^{k}(n_i-1)$. When $ p $ is fixed and $ min_{1<i<k} \,\,n_{i}\rightarrow \infty $, the asymptotic null distribution of $\varphi M$ test is chi-squared $ \chi^{2} $ with $ df=(k-1)p(p+1)/2 $ degrees of freedom, where
\begin{align*}
\varphi = 1- \frac{2p^{2}+3p-1}{6(p+1)(k-1)}\Big(\sum_{i=1}^{k} \frac{1}{n_{i}-1}-\frac{1}{n-k}\Big).
\end{align*}
The Box's $ M $ test in (6) represents a ratio of the pooled determinant $ |\textit{\textbf{S}}| $ to the geometric mean of the determinants $ |\textit{\textbf{S}}_{i}|, \, (i=1, ...,k) $, this test is extremely sensitive to departures from normality. On the other hand, the $M$ test is valid only for $n_i > p$ for $1\leq i\leq k$ since the sample covariance matrix $\textit{\textbf{S}}_{i}$ is singular ( $|\textit{\textbf{S}}_{i}|=0$ )  with positive probability when $p\geq n_i$ ( Dykstra, 1970). But  $\textit{\textbf{1}}\textit{\textbf{S}}_{i}\textit{\textbf{1}}^{T}$ can be positive almost surely when  $p$ is much bigger than $n_i$. Another disadvantage of the  $M$ test is that
it is not easy to calculate the determinant  $|\textit{\textbf{S}}|$ when $p$ is large. In order to overcome these two shortcomings, we  propose a modified $M$ test, $L_k$, by replacing
$|\textit{\textbf{S}}_{i}|$ in $ M $ test (6) with $ \textit{\textbf{1}}\textit{\textbf{S}}_{i}\textit{\textbf{1}}^{T}$, that is,
\begin{align}\label{8}
 L_{k}:= (n-k) log \hat{\textit{\textbf{V}}}_{p} - \sum_{i=1}^{k}(n_{i}-1)log \hat{\textit{\textbf{S}}}_{i}
\end{align}
where
\begin{align}\label{11}
\hat{\textit{\textbf{S}}}_{i}= \textit{\textbf{1}}\textit{\textbf{S}}_{i}\textit{\textbf{1}}^{T},\,\,\,\,\,\,\,\,\,\,\,\,\,\,\,  \hat{\textit{\textbf{V}}_{p}}=\frac{1}{n-k} \sum_{i=1}^{k} (n_{i}-1)\hat{\textit{\textbf{S}}_{i}},\,\,\,  i=1,...,k.
\end{align}
It can be seen that the modified $M$ test, $L_k$, can not only be defined but also can be easily calculated when $p$ is large.

\section{Main results}

The asymptotic properties of $\textit{\textbf{1}}\textit{\textbf{S}}_{i}\textit{\textbf{1}}^{T}$ and $L_k$ will be given in this section.

\par
\textbf{Theorem 3.1.}\textit{ For each $i$, $1\leq i\leq k$, let $ \{X^{(i)}_{jl}, 1\leq j\leq p\}, 1\leq l\leq n_i$, be $n_i$ i.i.d. G-stationary processes with $ E(X^{(i)}_{jl})=0$ and $ E(|X_{jl}^{(i)}|^r)<\infty $ for $ r>2,   1\leq l \leq n_i $, $1\leq j\leq p$ and $1\leq i\leq k$. Let $ \textbf{1}\Sigma_i \textbf{1}^{T}>0 $ and $\sigma^2_i=\sum_{j\in Z}\gamma^{(i)}(j)$ for $1\leq i\leq k$. If $p\geq c(n_i)^{r/(r-2)}$ for some positive constant $c$, then
\begin{align}\tag{9}\label{9}
\frac{(n_{i}-1) \textbf{1}\,\textit{\textbf{S}}_{i}\,\textbf{1}^{T}}{p\sigma^2_i} \to  \chi_i, \, a.s.,\,\,\,\,\,\, \chi_i \sim \chi_{(n_{i}-1)}^{2}
\end{align}
as $p\to \infty$ for $1\leq i\leq k.$}

\par
\textbf{Remark 3.1.} From (3) it follows that $\textbf{1}\Sigma_i \textbf{1}^{T}/(p\sigma^2_i) \to 1$  as $p\to \infty$, that is, $(n_{i}-1) \textbf{1}\,\textit{\textbf{S}}_{i}\,\textbf{1}^{T}/\textbf{1}\Sigma_i \textbf{1}^{T} \to  \chi_i$.  Gupta and Nagar (2000, Theorem 3.3.11) has proved the similar result of Theorem 3.1 for $ p< n_i $ by using the method of characteristic function assuming the normal distribution of the data. Here, we prove Theorem 3.1. without assuming the normality and permitting  $ p $ larger than $n_i$, and therefore (9) holds for any data distribution satisfying the conditions of Theorem 3.1.

Let $n_{min}=\min\{n_1, n_2,...,n_k\}$ and $n_{max}=\max\{n_1, n_2,...,n_k\}$. The following theorem shows that the asymptotic distribution of the test statistic $L_{k}$ in (7) is $\chi^{2}_{k-1}$.

\par
\textbf{Theorem 3.2.} \textit{ Let  the conditions of Theorem 2.1 hold and $p \geq c(n_{max})^{3r/(2-r)}$. Under the original hypothesis $ H_{0}$ in (1) with  $\Sigma >0$, we have that the test statistic $L_{k}$ in (7)  converges in distribution to $\chi^{2}_{k-1}$ as $n_{min} \to \infty$.}

Next we give an example to illustrate how to use Theorem 3.1. to test the null hypothesis in (1).

\textbf{Example 1}. Let $k=3$ and $n_1=n_2=n_3=100$. Without loss of generality, we assume that $p$ is even.  Let $\alpha =0.05$. Take two positive numbers $\chi_1$ and $\chi_2$ such that $\chi_1 <\chi_2$ and
\begin{eqnarray*}
\textbf{P}( \chi^{2}_{k-1} < \chi_1)=0.025,\,\,\,\,\,\,\,\,\,\,\, \textbf{P}( \chi^{2}_{k-1} > \chi_2)=0.025.
\end{eqnarray*}
This means that $\textbf{P}( \chi_1\leq  \chi^{2}_{k-1} \leq  \chi_2)=1-0.05$, that is, $\textbf{R}_k=\{ [0, \, \chi_1)\cup ( \chi_2,\, \infty)\}$ is a domain of rejecting the null hypothesis $H_0$,
here, the procedures of $L_{k}$ hold for any $ \textit{\textbf{1}}=(1_1,..., 1_p) $, $ \textit{\textbf{1}} $ is $ (1 \times p) $ vector, .

For high dimensional observation vector $(X^{(i)}_1, ..., X^{(i)}_p), 1\leq i\leq k$, researchers sometimes not only want to infer whether their covariance matrices, $\Sigma_1, ..., \Sigma_k,$ are  equal, but also want to know whether sub covariance matrices of relative low dimensional observations vectors,  $ (X^{(i)}_1, ..., X^{(i)}_{p_1})$, $ (X^{(i)}_{p_1+1}, ..., X^{(i)}_{p_1+p_2})$, ..., $ (X^{(i)}_{p_1+...p_{l-1}+1}, ..., X^{(i)}_{p})$, are equal, where  $p=\sum_{i=1}^lp_i$, all $p_i\geq 2$ and $l\geq 2$.

Without loss of generality,  we assume that $p \geq n_{min}=\min\{n_i, 1\leq i\leq k\}$. Let $p_j=j (n_{min}-1)\wedge p$ for $j\geq 0$ and take $m\geq 2$ such that $p_{m-1}< p \leq  p_m$. By the definition of $p_j$, we know that $p_m=p$.  Now, we divide $\Sigma_i$ and $\textit{\textbf{S}}_{i}$ into the following two forms respectively
\begin{eqnarray*}
\Sigma_i=\left(\begin{array}{cccccccc} \Sigma^{(i)}_{p_1\times p_1} & \Sigma^{(i)}_{p_1\times p_2}&  ...&...&\\
\Sigma^{(i)}_{p_2\times p_1} & \Sigma^{(i)}_{p_2\times p_2} & ... &...&\\
...&...&...&...&...\\
...&...&...&...&\Sigma^{(i)}_{p_{m-1}\times p}\\ &...& &\Sigma^{(i)}_{p\times p_{m-1}}& \Sigma^{(i)}_{p\times p}
\end{array}\right)
\end{eqnarray*}
and
\begin{eqnarray*}
\textit{\textbf{S}}_{i}=\left(\begin{array}{cccccccc} S^{(i)}_{p_1\times p_1} & S^{(i)}_{p_1\times p_2}&  ...&...&\\
S^{(i)}_{p_2\times p_1} & S^{(i)}_{p_2\times p_2} & ... &...&\\
...&...&...&...&...\\
...&...&...&...&S^{(i)}_{p_{m-1}\times p}\\ &...& &S^{(i)}_{p\times p_{m-1}}& S^{(i)}_{p\times p}
\end{array}\right),
\end{eqnarray*}
for $1\leq i\leq k$, where $\Sigma^{(i)}_{p_j\times p_j}$ and  $S^{(i)}_{p_j\times p_j}$ are all $(p_{j}-p_{j-1})\times (p_{j}-p_{j-1})$ sub-matrices for $1\leq i\leq k$ and $1\leq j\leq m-1$, and $\Sigma^{(i)}_{p\times p}$ and  $S^{(i)}_{p\times p}$ are $(p-p_{m-1})\times (p-p_{m-1})$ sub-matrices for $1\leq i\leq k$.

Let $\textbf{y}_{p_1}=(1_1,..., 1_{p_1}), \textbf{y}_{p_2}=(1_{p_1+1},..., 1_{p_2}),..., \textbf{y}_{p_m}=(1_{p_{m-1}+1},..., 1_{p})$ and
\begin{align}\tag{10}\label{10}
\textbf{y}_1=(\textbf{y}_{p_1}, 0, ...,0), \textbf{y}_2=(0,...,0, \textbf{y}_{p_2},0, ...,0), ..., \textbf{y}_m=(0,...,0,\textbf{y}_{p_m}),
\end{align}
be $m$ $p$-dimensional vectors. It follows that
\begin{eqnarray*}
\textbf{y}_j\,\Sigma_i\,\textbf{y}_j^{T}=\textbf{y}_{p_j}\Sigma^{(i)}_{p_j\times p_j}\textbf{y}_{p_j}^T,\,\,\,\,\,\,\,\,\textbf{y}_j\,\textit{\textbf{S}}_{i}\,\textbf{y}_j^{T}=\textbf{y}_{p_j}S^{(i)}_{p_j\times p_j}\textbf{y}_{p_j}^T
\end{eqnarray*}
for $1\leq i\leq k, \, 1\leq j\leq m$.  Note that $p_j-p_{j-1}=n_{min}-1$ for $1\leq j\leq m-1$ and $p-p_{m-1}\leq n_{min}-1$, that is, $\textbf{y}_{p_j}, 1\leq j\leq m-1,$ are $(n_{min}-1)$-dimensional vectors and $\textbf{y}_{p_m}$ is a  $(p-p_{m-1})$-dimensional vector.  If $\textbf{y}_j \neq 0$, or equivalently, $\textbf{y}_{p_j}\neq 0$, then $\textbf{y}_j\,\textit{\textbf{S}}_{i}\,\textbf{y}_j^{T}=\textbf{y}_{p_j}S^{(i)}_{p_j\times p_j}\textbf{y}_{p_j}^T>0$  when $\Sigma^{(i)}_{p_j\times p_j}>0$ for $1\leq i\leq k, \, 1\leq j\leq m$.  Now we present $m$ test statistics as follows:
\begin{align}\tag{11}\label{11}
 L_{kj} &:= (n-k) log \hat{\textit{\textbf{V}}}_{pj} - \sum_{i=1}^{k}(n_{i}-1)log \hat{\textit{\textbf{S}}}_{ij}
\end{align}
for $1\leq j\leq m$, where
\begin{align*}
\hat{\textit{\textbf{S}}}_{ij}=\textbf{y}_j\,\textit{\textbf{S}}_{i}\,\textbf{y}_j^{T}=\textbf{y}_{p_j}S^{(i)}_{p_j\times p_j}\textbf{y}_{p_j}^T, \,\,\,\,\,\,\, \,\,\, \hat{\textit{\textbf{V}}_{pj}}=\frac{1}{n-k} \sum_{i=1}^{k} (n_{i}-1)\hat{\textit{\textbf{S}}_{ij}}
\end{align*}
for $1\leq i\leq k$ and $1\leq j\leq m$.

Hence, by extending Theorem 3.2 we can get the following corollary.

\par \textbf{Corollary 3.1.}  \textit{ Let  $p \geq c(n_{max})^{3r/r-2}$. Under the original hypothesis $ H_{0}$ in (1) with  $\Sigma >0$, for any $p$-dimensional vector $\textbf{y}_{j} \neq 0$ defined in (10), we have that for every $j$ ($1\leq j\leq m$), the test statistic $L_{kj}$ in (11) converges in distribution to $\chi^{2}_{k-1}$ as $n_{min} \to \infty$.}

\par \textbf{Remark 3.2.} Here, we divide the original hypothesis $ H_{0}$ and the alternative hypothesis
$H_1$ into $m$ original hypotheses and $m$ alternative hypotheses in the following
\begin{align*}
 H_{0j}: \Sigma^{(1)}_{p_j\times p_j}=...= \Sigma^{(k)}_{p_j\times p_j}  \,\,\, vs\,\,\,  H_{1j}:  \Sigma^{(i)}_{p_j\times p_j} \neq \Sigma^{(l)}_{p_j\times p_j} \text{\,\, at least one pair $(i, l)$, $i \neq l$}
\end{align*}
for $1\leq j\leq m$. Thus, we can use the test statistic $L_{kj}$ to do the  original hypothesis $H_{0j}$ and the alternative hypothesis $H_{1j}$ for $1\leq j\leq m$.  Moreover, though the above division may lose some information, for example, $\Sigma^{(k)}_{p_j\times p_l}, \, j\neq l$, one of its advantages is that it can help us to test which part ( $H_{0j}$ ) of the original hypothesis $ H_{0}: \Sigma_{1}=\Sigma_{2}=...=\Sigma_{k}=\Sigma $ is inconsistent.

\par \textbf{Remark 3.3.} If necessary, one can divide the original hypothesis $ H_{0}$ and the alternative hypothesis $H_{1}$ into more or less sub hypotheses, $H_{0j}$ and $H_{1j}$, $1 \leq j\leq l$.

The following example will discuss on how to use Corollary 3.1. for testing the null hypothesis.

\textbf{Example 2}. Let $k=3$, $n_1=n_2=n_3=101$ and $p=350$.  Let $\textbf{1}_a=(1, 1, ..., 1)$ and $\textbf{0}_a=(0, 0, ..., 0)$ be two $100$-dimensional vectors and $\textbf{1}_b=(1, 1, ..., 1)$ and $\textbf{0}_b=(0, 0, ..., 0)$ denote two $50$-dimensional vectors. Take four $350$-dimensional vectors $\textbf{y}_j, 1\leq j\leq 5,$ defined in the following
\begin{eqnarray*}
&&\textbf{y}_1 =(\textbf{1}_a, \textbf{0}_a,  \textbf{0}_a, \textbf{0}_b), \,\,\,\,\,\,\,\,\,\,\textbf{y}_2=(\textbf{0}_a, \textbf{1}_a, \textbf{0}_a, \textbf{0}_b)\\
&&\textbf{y}_3 =(\textbf{0}_a, \textbf{0}_a,  \textbf{1}_a, \textbf{0}_b), \,\,\,\,\,\,\,\,\,\,\textbf{y}_4=(\textbf{0}_a, \textbf{0}_a, \textbf{0}_a, \textbf{1}_b).
\end{eqnarray*}
Divide  $\textbf{S}_i$ into four sub sample covariance matrices: $ S^{(i)}_{p_j\times p_j},  j=1, 2, 3, 4$  for $i=1, 2, 3$, where $p_1=100, p_2=200, p_3=300$ and $p_4=350$. By (7) we can get four test statistics $ L_{kj}, 1\leq j\leq 4$. Thus, we will reject the null hypothesis $H_0$ if at least one of $ L_{kj}, 1\leq j\leq 4$ belong to the rejecting region $\textbf{R}_k$.
If $L_{k4}\in \textbf{R}_k$ and $L_{kj} \notin \textbf{R}_k$ for $1\leq j\leq 3$, this means that the three covariances, $\Sigma^{(1)}_{p_4\times p_4}$, $\Sigma^{(2)}_{p_4\times p_4}$ and $ \Sigma^{(3)}_{p_4\times p_4} $, corresponding to $(\textbf{X}^{(1)}_{301}, ..., \textbf{X}^{(1)}_{350})$, $(\textbf{X}^{(2)}_{301}, ..., \textbf{X}^{(2)}_{350})$ and $(\textbf{X}^{(1)}_{301}, ..., \textbf{X}^{(1)}_{350})$,  are not equal.

\par \textbf{Remark 3.4.} In practice, $ \rho L_{k} $ converges to $ \chi^{2}_{k-1} $, where $ \rho=\rho_{n} $ is the scale factor given in Bartlett (1937) as follows:
\begin{align}\tag{12}\label{12}
\rho=\rho_{n}=\frac{1}{C}, \,\,\, where \,\,\, C= 1+\frac{1}{3(k-1)}\Big(\sum_{i=1}^{k} \frac{1}{n_{i}-1}-\frac{1}{n-k}\Big)
\end{align}
and  $ \rho \rightarrow 1 $ for fixed $ k $ as  $ n_{min}\rightarrow \infty $ . Here, the use of $ \rho $ tends to over-correct the exaggerated significance levels in small samples. The scale factor $ \rho $ also greatly increases the approximation of significance level and uses as a gauge of its convergence. The empirical sizes of $ \rho L_{k} $ in Table 1 are corresponding to the correct significance levels of $ \chi^{2}_{k-1} $ for several cases of $ n_{i} $, $ p $ and $ k $ with the nominal level 0.05.

\section{Numerical simulations}
In this section, we use $ \rho L_{k} $ to test the hypothesis (1). Here, $ \rho $ was chosen as in (12) so that the convergence of the significance levels by the empirical distributions of $ \rho L_{k} $ test are close to $ \chi^{2}_{k-1} $. All simulation results in this paper are obtained using $ 10^{3} $ repetitions with the nominal level $ (0.05) $.
To demonstrate how $ \rho L_{k} $ performs, we conducted the test proposed by Ahmad (2017) ($ T_{g} $ test) in non-normal high dimensional data. We compare the power and empirical size of the $ \rho L_{k} $ with that of the   $ T_{g} $, where $ T_g= (g-1)\sum_{i=1}^g T_i -2 \sum_{i=1}^g \sum_{j=1}^g T_{ij} $, see Ahmed (2017).

\paragraph{•} \textbf{The empirical sizes and attained power }

Based on different distributions structures under the null hypothesis, our simulations
are designed in Two scenarios: (i) In dependent stationary process of Gaussian AR(1) model, to evaluate the empirical sizes and power of the $ \rho L_{k} $  tests, we present the results
of the equality test between three covariance matrices $ (k=3) $. For this purpose,
we generated three multivariate normal random vectors $  \textbf{X}^{(i)}_{jl_G},\,\, , 1\leq j\leq p,  1\leq i\leq 3, 1\leq l\leq n_i  $,
with $ \overline{\textbf{X}}^{(i)}_{jl_G} =\mu^{(i)}_{p}= 0_p $ and  $ cov(\textbf{X}^{(i)}_{jl_G})=\Sigma_{i} =\, \Sigma^{(i)}_{p} = \, (\sigma_{i_{\cdot}j_{\cdot}})_{p\times p}
 =(-1)^{(i_{\cdot}+j_{\cdot})}(0.2 \times (J+2))^{\vert i_{\cdot}-j_{\cdot}\vert}\, = \Omega_J, \, J=0,1,2 $. We considered the following hypothesis testing setup:

$ H_0: \Sigma_{1}=\Sigma_{2}=\Sigma_{3}=\Sigma = \Omega_0 $, $ H_1: \Sigma_{1}=  \Sigma $, $\Sigma_{2}= \Omega_1 $ and $  \Sigma_{3}= \Omega_2 $.\\
(ii) Here, the dependent stationary processes of multivariate uniform distribution $ \textbf{X}^{(i)}_{jl_U} $ and multivariate exponential distribution $ \textbf{X}^{(i)}_{jl_E} $ were generated by AR(1) model. For testing the null hypotheses we adapted the foregoing setting of $ H_0: \Sigma_{1}=\Sigma_{2}=\Sigma_{3}=\Sigma = \Omega_0 $ in (i) above for both the uniform and exponential processes $ \textbf{X}^{(i)}_{jl_U} $ and $ \textbf{X}^{(i)}_{jl_E} $, respectively. Moreover,
for power, the alternative hypotheses were constructed at least one of the three distributions that follows a different covariance structure than the others. The first two structures of covariance matrices are defined as Compound Symmetry (CS):= $\Sigma_{1}= K \textbf{I}+\phi \textbf{J} $ and Autoregressive of order 1 (AR(1)):=  $ \Sigma_{2} =(-1)^{(i_{\cdot}+j_{\cdot})}(0.4))^{\vert i_{\cdot}-j_{\cdot}\vert}\ $, respectively, with $ \textbf{I} $ as identity matrix, $ \textbf{J} $ as the matrix of ones and $ K=1 $, $ \phi=0.70 $ as constants, However, the third covariance matrix $ \Sigma_{3} $ was simulated with data being drawn from the p-dimensional centralized $ \Gamma (1,1) $.

We computed the empirical size under  $ H_0 $ and the attained power under $ H_1 $,  Tables 1 and 2 show the results of $ \rho L _{k} $ and $T_g$, respectively.
From Table 1 we find that the $ \rho L_{k} $ test has good performance for the three distributions.
The empirical sizes of $ \rho L_{k} $ are uniformly close to the nominal level $ 0.05 $.
We observe that the power of $ \rho L _{k} $ test uniformly approaches to 1 for small, moderate and large sample sizes $ n_{i} $. According to Table 2, we observe that the performance
of $ \rho L_{k} $ test is better than that of $ T_{g} $ test in all cases.
The reason why the performance of the proposed test is better should be that the test statistic
$ \rho L_{k} $ is based on the $k$ estimators $ \textbf{\textit{1}}\textit{\textbf{S}}_{i}\textbf{\textit{1}}^{T}$, $ 1\leq i\leq k. $

\begin{table}[H]
\begin{center}
\caption{ The empirical size and power of $ \rho L_{k} $ test in the normal, exponential and uniform  distributions}
\label{table:table-1}
\begin{tabular}{|c|c|ccc|ccc|}\hline
	& &  \multicolumn{3}{c|}{Empirical size}  &\multicolumn{3}{c|}{Pwoer} \\\hline
  $ p $& $ n_{1} = n_{2}=n_{3} $  & Normal & Exp & Uniform  & Normal & Exp & Uniform \\\hline
 20 & 10 & 0.049  & 0.054 & 0.057 &     0.981 & 0.988  &  0.984   \\
    & 20 & 0.049  & 0.056 & 0.056 &     0.989 & 0.999  &  0.993 \\
    & 50 & 0.051  & 0.047 & 0.061 &     1.000 & 1.000  &  1.000  \\
    & 100& 0.050  & 0.054 & 0.067 &     1.000 & 1.000  &  1.000 \\
    \hline

 50 & 10 & 0.054  & 0.052 & 0.052 &     0.992 & 0.999  &  0.990\\
    & 20 & 0.046  & 0.050 & 0,049 &     0.991  & 1.000  &  1.000\\
    & 50 & 0.045  & 0.047 & 0.056 &     1.000  & 1.000  &  1.000 \\
    & 100& 0.053  & 0.050 & 0.051 &     1.000  & 1.000  &  1.000 \\
    \hline

 100& 10  & 0.050 & 0.050 & 0.051 &     0.989  & 1.000  &  1.000\\
    & 20  & 0.048 & 0.042 & 0.052 &     1.000  & 1.000  &  1.000 \\
    & 50  & 0.053 & 0.054 & 0.046 &     1.000  & 1.000  &  1.000 \\
    & 100 & 0.052 & 0.050 & 0.052 &     1.000  & 1.000  &  1.000\\
    \hline

 200& 10 & 0.047  & 0.046 & 0.041 &     0.987  & 1.000  &  1.000\\
    & 20 & 0.046  & 0.050 & 0.055 &     1.000  & 1.000  &  1.000\\
    & 50 & 0.050  & 0.047 & 0.053 &     1.000  & 1.000  &  1.000 \\
    & 100& 0.051  & 0.052 & 0.050 &     1.000  & 1.000  &  1.000\\
    \hline
 300& 10 & 0.053  & 0.051 & 0.054 &     0.995  & 1.000  &  1.000\\
    & 20 & 0.046  & 0.050 & 0,048 &     1.000  & 1.000  &  1.000\\
    & 50 & 0.050  & 0.047 & 0.051 &     1.000  & 1.000  &  1.000\\
    & 100& 0.051  & 0.051 & 0.053 &     1.000  & 1.000  &  1.000\\
    \hline

\end{tabular}
\end{center}
\end{table}

\begin{table}[H]
\begin{center}
\caption{ The empirical size and power of $ T_{g} $ test in the normal, exponential and uniform  distributions}
\label{table:table-1}
\begin{tabular}{|c|c|ccc|ccc|}\hline
	& &  \multicolumn{3}{c|}{Empirical size}  &\multicolumn{3}{c|}{Pwoer} \\\hline
  $ p $& $ n_{1} = n_{2}=n_{3} $  & Normal & Exp & Uniform  & Normal & Exp & Uniform \\\hline
 20 & 10 & 0.001  & 0.005 & 0.010 &     0.482 & 0.201  &  0.320   \\
    & 20 & 0.019  & 0.006 & 0.012 &     0.518 & 0.190  &  0.300 \\
    & 50 & 0.021  & 0.011 & 0.015 &     0.540 & 0.234  &  0.341  \\
    & 100& 0.023  & 0.014 & 0.017 &     0.618 & 0.261  &  0.365 \\
    \hline

 50 & 10 & 0.014  & 0.010 & 0.012 &     0.452  & 0.251  &  0.380\\
    & 20 & 0.024  & 0.015 & 0,016 &     0.672  & 0.351  &  0.480\\
    & 50 & 0.023  & 0.020 & 0.019 &     0.751  & 0.391  &  0.590 \\
    & 100& 0.025  & 0.022 & 0.023 &     0.800  & 0.421  &  0.680 \\
    \hline

 100& 10  & 0.016 & 0.013 & 0.015 &     0.513  & 0.301  &  0.490\\
    & 20  & 0.024 & 0.014 & 0.019 &     0.712  & 0.362  &  0.602 \\
    & 50  & 0.026 & 0.017 & 0.021 &     0.821  & 0.401  &  0.680 \\
    & 100 & 0.029 & 0.020 & 0.022 &     0.850  & 0.510  &  0.750\\
    \hline

 200& 10 & 0.017  & 0.014 & 0.019 &     0.590  & 0.452  &  0.522\\
    & 20 & 0.016  & 0.017 & 0.020 &     0.763  & 0.562  &  0.682\\
    & 50 & 0.028  & 0.020 & 0.023 &     0.862  & 0.572  &  0.702 \\
    & 100& 0.030  & 0.022 & 0.025 &     0.923  & 0.662  &  0.792\\
    \hline
 300& 10 & 0.021  & 0.019 & 0.022 &     0.640  & 0.500  &  0.541\\
    & 20 & 0.027  & 0.021 & 0.024 &     0.792  & 0.782  &  0.727\\
    & 50 & 0.027  & 0.020 & 0.025 &     0.831  & 0.702  &  0.771\\
    & 100& 0.030  & 0.024 & 0.029 &     0.913  & 0.762  &  0.802\\
    \hline

\end{tabular}
\end{center}
\end{table}

\section{Experimental study}
\paragraph{•} To demonstrate the performance of the proposed  $ \rho L_k $ test we use
the time-series of LSVT Voice rehabilitation dataset, provided by https://archive.ics.uci.edu.
This dataset is used to study the extracting of clinical information from speech signals  provided by LSVT Global, a company specializing in voice rehabilitation (see Tsanas, M.A. et al. 2014). The original study used 309 Attributes (variables) to characterize 126 speech signals (Instances).

\subsection*{Results and discussion }
We show the results of the LSVT voice rehabilitation data set in Table 3.
We chose three sample sizes $ n_{i}=10, n_{i}=20 $ and $ n_{i}=40 $,
each sample size was chosen by randomly selecting instances.
Table 5 show the statistics measure and
$ P_{value} $ of $ T_{g}$ and $ \rho L_{k} $ tests.
We observe in Table 3 that the $ \rho L_{k} $ test rejects the null hypothesis $ H_0 $, that is, $ P_{value} <0.05 $. The results of  $ \rho L_k $ may consistent with the results of $ T_{g}$, however, for $ n_i=10 $ the $ P_{value}$ of $ \rho L_{k} $ test equals $0.0391 $, that is, $ \rho L_{k} $ appears more conservative than the $ T_{g}$ to reject $H_0$ when $H_0$ is correct.

\begin{table}[H]
\begin{center}
\caption{ The statistic measure and $ P_{value} $  of $ T_{g}$ and $ \rho L_{k} $ tests, of LSVT voice rehabilitation dataset}
\label{table:table-3}
\begin{tabular}{|c|cc|cc|}\hline
  & \multicolumn{2}{c|}{$ T_{g} $}  &\multicolumn{2}{c|}{$ \rho L_{k} $}   \\\hline
 $ n_{1} = n_{2}=n_{3} $   &  $ (z)_T $  & $ P_{value} $  & $ (\chi^{2})_{\rho L} $  & $ P_{value} $    \\\hline
10 &  4.2005  & 0.0000   & 3.0833   &0.0391   \\
20 &  9.4051  & 0.0000   & 15.1066  &0.0005     \\
40 &  13.0521 & 0.0000   & 39.7810  & 0.0000  \\
    \hline
\end{tabular}
\end{center}
\end{table}

\section{Concluding remarks}

This article presents a modified test, $L_{k} $, to test the null hypothesis in (1) for high-dimensional settings with ignoring the assumption of normality. We prove that the $L_{k} $ converges in distribution to $\chi_{k-1}^{2}$. Here, we don't assume the usual conditions about the relationship between the traces of powers of covariance matrix and $p$,
but we assume that the relationship between $p$ and $n_{max}$ satisfies  $p\geq c(n_{max})^{3r/(r-2)} $ for a positive constant $c$.
Hence, the $L_k$ test accommodates the large $p$, small $n_{max}$ situations. The most significant two advantages of the modified test, $L_{k} $, are that (1) it is easy to calculate; (2) it can effectively test the equality of several high-dimensional covariance matrices when $p$ is much bigger than $n_{max}$.  The simulation results in Tables 1-2 demonstrate that the $\rho L_{k}$ test has perfect performance and its power tends to 1.

\newpage
\section{Appendix}

\paragraph{•} \textbf{A.1. Proof of Theorem 3.1.}
Let $\mathbb{S}^{(i)}_{pl}=\sum_{j=1}^pX^{(i)}_{jl}$ and $\overline{\mathbb{S}^{(i)}}_{p}=\sum_{l=1}^{n_i}\mathbb{S}^{(i)}_{pl} /n_i$. By the strong invariance principle in (2), $\mathbb{S}^{(i)}_{pl}$ can be written as
\begin{align}\tag{A.1}\label{A.1}
\mathbb{S}^{(i)}_{pl}=\sigma_i\mathbb{B}^{(i)}_l(p)+\epsilon^{(i)}_{pl}, \,\,\,\, a.s.,
\end{align}
for $1\leq l\leq n_i, 1\leq i\leq k$, and therefore,
\begin{align}\tag{A.2}\label{A.2}
\mathbb{S}^{(i)}_{pl}-\overline{\mathbb{S}^{(i)}}_{p}=\sigma_i\Big(\mathbb{B}^{(i)}_l(p)-\overline{\mathbb{B}^{(i)}}(p)\Big)+\epsilon^{(i)}_{pl}-\overline{\epsilon^{(i)}}_{p}, \,\,\,\, a.s.,
\end{align}
where $\overline{\mathbb{B}^{(i)}}(p)=\sum_{l=1}^{n_i}\mathbb{B}^{(i)}_l(p)/n_i$, for each $i$, $\mathbb{B}^{(i)}_l(t), 1\leq l\leq n_i,$ are mutually independent $n_i$ standard Brownian motions, and the random variables $\epsilon^{(i)}_{pl}$ satisfies that $\epsilon^{(i)}_{pl}/p^{1/r} \to 0, \,\, a.s.$ as $p\to \infty$ for all $1\leq l\leq n_i, 1\leq i\leq k$. Note that
\begin{align}\tag{A.3}\label{A.3}
\chi_i:=\frac{\sum_{l=1}^{n_i}(\mathbb{B}^{(i)}_l(p)-\overline{\mathbb{B}^{(i)}}(p))^2}{p}\,\, \sim \,\, \chi^{2}_{(n_i-1)}.
\end{align}
It means that $\chi_i$ does not depend on $p$ for all $1\leq l\leq n_i$,  $n_i>1$ and $1\leq i\leq k$. From (5), (A.1) and (A.2), it follows that
\begin{align}\tag{A.4}\label{A.4}
\frac{(n_{i}-1) \textbf{1}\,\textit{\textbf{S}}_{i}\,\textbf{1}^{T}}{p\sigma^2_i} =\frac{1}{p\sigma^{2}}\sum_{l=1}^{n_{i}}(\mathbb{S}^{(i)}_{pl}-\overline{\mathbb{S}^{(i)}}_{p})^2=\chi_i+\xi_i(p)
\end{align}
for large $p$, $1\leq l\leq n_i$,  $n_i>1$ and $1\leq i\leq k$, where
\begin{align}\tag{A.5}\label{A.5}
\xi_i(p)&=2\frac{\sum_{l=1}^{n_i}(\mathbb{B}^{(i)}_l(p)-\overline{\mathbb{B}^{(i)}}(p))}{\sqrt{p}}\frac{(\epsilon^{(i)}_{pl}-\overline{\epsilon^{(i)}}_{p})}{\sqrt{p}}+\frac{1}{p\sigma^2_i}\sum_{l=1}^{n_i}(\epsilon^{(i)}_{pl}-\overline{\epsilon^{(i)}}_{p})^2\\
&=\frac{\sqrt{n_i}p^{1/r}}{\sqrt{p}}\Big(\frac{\sum_{l=1}^{n_i}(\mathbb{B}^{(i)}_l(p)-\overline{\mathbb{B}^{(i)}}(p))}{\sqrt{n_ip}}\frac{(\epsilon^{(i)}_{pl}-\overline{\epsilon^{(i)}}_{p})}{p^{1/r}}+\frac{\sqrt{n_i}p^{1/r}}{\sqrt{p}}\frac{1}{n_i}\sum_{l=1}^{n_i}\frac{(\epsilon^{(i)}_{pl}-\overline{\epsilon^{(i)}}_{p})^2}{p^{2/r}}\Big).\nonumber
\end{align}
Since $(\mathbb{B}^{(i)}_l(p)-\overline{\mathbb{B}^{(i)}}(p))/\sqrt{n_ip} \, \sim N(0, (n_i-1)/n^2_i)$, $(\epsilon^{(i)}_{pl}-\overline{\epsilon^{(i)}}_{p})/p^{1/r} \to 0, \, a.s.$ as $p\to \infty$ and $\sqrt{n_i}p^{1/r}/\sqrt{p}$ is bounded, or  $p\geq c(n_i)^{r/(r-2)}$ for some positive constant $c$, it follows that $\xi_i(p) \to 0, \, a.s.,$ as $p\to \infty$ for $1\leq i\leq k$. Thus, by (A.4) and (A.5) we have
\begin{eqnarray*}
\frac{(n_{i}-1) \textbf{1}\,\textit{\textbf{S}}_{i}\,\textbf{1}^{T}}{p\sigma^2_i}  \rightarrow \, \chi_i, \, a.s. \,\,\,\,\,\,\,\, \chi_i \sim \chi^{2}_{n_i-1}  
\end{eqnarray*}
as $p\to \infty$ for $1\leq i\leq k$.

\paragraph{•} \textbf{A.2. Proof of Theorem 3.2.}  Under the original hypothesis $ H_{0}$ in (1) with  $\Sigma >0$, we have $\sigma^2_1=\sigma^2_2=...=\sigma^2_k=\sigma^2$. Let $ a_{i}= (n_{i}-1)/(n-k) $ and $ N=(n-k)/2$. By (A.4) $ L_{k}$ can be written as $L_{k}=-2 Nlog J_{n,k} = -2log J^{^{N}}_{n,k} $ for $ 1 \leq i \leq k $, where
\begin{eqnarray*}
 J_{n,k} &=& \frac{\prod_{i=1}^{k}\Big[\hat{\textit{\textbf{S}}}_{i}\Big]^{\frac{n_{i}-1}{n-k}}}{\hat{\textit{\textbf{V}}}_{p}} = \frac{(n-k)\prod_{i=1}^{k}\Big[\hat{\textit{\textbf{S}}}_{i}\Big]^{a_{i}}\,\,\, \frac{1}{p\sigma^2}}{\sum_{i=1}^{k}(n_{i}-1)\frac{\hat{\textit{\textbf{S}}}_{i}}{p\sigma^2}} = \frac{(n-k)\prod_{i=1}^{k}\Big[\hat{\hat{\textit{\textbf{S}}}}_{i}\Big]^{a_{i}}}{\sum_{i=1}^{k}(n_{i}-1)\hat{\hat{\textit{\textbf{S}}}}_{i}} \\
 &=&  \frac{\prod_{i=1}^{k}\Big[a_{i}\Big]^{-a_{i}} \prod_{i=1}^{k} \Big[(n_{i}-1) \hat{\hat{\textit{\textbf{S}}}}_{i}\Big]^{a_{i}}}{\sum_{i=1}^{k}(n_{i}-1)\hat{\hat{\textit{\textbf{S}}}}_{i}} = C_{n,k}^{-1}\prod_{i=1}^{k}\Big[\frac{\chi_{i}Y_i}{\sum_{i=1}^{k} \chi_{i}}\Big]^{a_{i}}= C_{n,k}^{-1}\prod_{i=1}^{k}(Z_iY_i)^{a_{i}}
\end{eqnarray*}
where $\hat{\hat{\textit{\textbf{S}}}}_{i}:=\hat{\textit{\textbf{S}}}_{i}/p\sigma^2$, 
\begin{eqnarray*}
&&C_{n,k}=\prod_{i=1}^{k} \Big[\frac{n_{i}-1}{n-k}\Big]^{\frac{n_{i}-1}{n-k}}, \,\,\,\,\,\,\,\,\, Y_i=\frac{1+\xi_i(p)/\chi_{i}}{1+\sum_{i=1}^k\xi_i(p)/\sum_{i=1}^k\chi_{i}},\,\,\,\, 1\leq i\leq k,\\
&&Z_{i}=\frac{\chi_{i}}{\chi_{1}+...+\chi_{k}}, \,\,\,\,\,\,\,\,\, 1\leqslant i \leqslant k-1, \,\,\,\,\,\,\, Z_{k}=1-\sum_{i=1}^{k-1} Z_{i}.
\end{eqnarray*}
Since $\chi_{i}  \sim \chi^{2}_{(n_{i}-1)} $  for $1\leq i\leq k$, it follows that $ (Z_{1},...,Z_{k-1}) $ is subject to the Dirichlet distribution with density function $f(z_1,...,z_{k-1})$ in the following
\begin{eqnarray*}
f(z_1,...,z_{k-1})=\frac{\Gamma(N)}{\prod_{i=1}^{k}\Gamma(Na_{i})}\,\prod _{i=1}^{k}z_{i}^{Na_{i}-1}.
\end{eqnarray*}
As in the lemma 2 in Chao and Glaser (1978), we can get that
\begin{align}\tag{A.6}\label{A.6}
\textbf{E}\Big(\Big[C_{nk}^{-1} \prod_{i=1}^{k} (Z_{i})^{a_{i}}\Big]^{N\gamma}\Big)=\Big(C_{nk}^{-1}\Big)^{N\gamma} \frac{\Gamma(N)}{\Gamma(N(1+\gamma))} \prod_{i=1}^{k} \frac{\Gamma(Na_{i}(1+\gamma))}{\Gamma(Na_{i})}.
\end{align}
for any constant $ \gamma $ satisfying $ \gamma > -1 $. By using Stirling approximation to the gamma function in (A.6), it follows that
\begin{align}\tag{A.7}\label{A.7}
\textbf{E}\Big(\Big[C_{nk}^{-1} \prod_{i=1}^{k} (Z_{i})^{a_{i}}\Big]^{N\gamma}\Big)\,\,\rightarrow \,\, \Big(1+\gamma \Big)^{\frac{-(k-1)}{2}} 
\end{align}
as $ n_{min}\rightarrow \infty$.

On the other hand, if $Y_i>1$, then
\begin{align}\tag{A.8}\label{A.8}
(Z_iY_i)^{n_i-1}&=[Z_i(Y_i-1)+Z_i]^{n_i-1}\nonumber\\
&=(Z_i)^{n_i-1}+(n_i-1)Z_i(Y_i-1)\sum_{j=1}^{n_i-1}C^{j-1}_{n_i-2}\frac{1}{j}[Z_i(Y_i-1)]^{j-1}(Z_i)^{n_i-2-(j-1)} \nonumber\\
&\leq (Z_i)^{n_i-1}+(n_i-1)Z_i(Y_i-1)[Z_i(Y_i-1)+Z_i]^{n_i-2}\nonumber\\
&=(Z_i)^{n_i-1}+(n_i-1)Z_i(Y_i-1)(Z_iY_i)^{n_i-2}, \,\,\,\,\, n_i \geq 2 \nonumber\\
(Z_iY_i)^{-\alpha}&\leq (Z_i)^{-\alpha}, \,\,\,\,\, \alpha >0 \nonumber
\end{align}
for $1\leq i\leq k$. Similarly, if $Y_i\leq 1$, then
\begin{align}\tag{A.9}\label{A.9}
(Z_iY_i)^{n_i-1}&=[Z_i-Z_i(1-Y_i)]^{n_i-1}\nonumber\\
&=(Z_i)^{n_i-1}-(n_i-1)Z_i(1-Y_i)\sum_{j=1}^{n_i-1}C^{j-1}_{n_i-2}\frac{1}{j}[-Z_i(1-Y_i)]^j(Z_i)^{n_i-2-(j-1)} \nonumber\\
&\geq (Z_i)^{n_i-1}-(n_i-1)Z_i(1-Y_i)[Z_i-Z_i(1-Y_i)]^{n_i-2}\nonumber\\
&\geq (Z_i)^{n_i-1}-(n_i-1)Z_i(1-Y_i)(Z_i)^{n_i-2}, \,\,\,\,\, n_i\geq 2\nonumber\\
(Z_iY_i)^{-\alpha}&\geq (Z_i)^{-\alpha}, \,\,\,\,\, \alpha >0 \nonumber
\end{align}
for $1\leq i\leq k$. 

Note that $0\leq Z_iY_i\leq 1$, $0<Z_i<1$, and therefore, $|Z_i(1-Y_i)|=|Z_i(Y_i-1)|\leq 1$. By (A.5) we know that
\begin{eqnarray*}
n_i|Y_i-1|=\frac{n_i|\xi_i(p)/\chi_{i}-\sum_{i=1}^k\xi_i(p)/\sum_{i=1}^k\chi_{i}|}{1+\sum_{i=1}^k\xi_i(p)/\sum_{i=1}^k\chi_{i}} \, \to \, 0, \,\,\, a.s.
\end{eqnarray*}
as $n_{min}\to \infty$ for $1\leq i\leq k$ when $n_{max}\sqrt{n_{max}}p^{1/r}/\sqrt{p}$ is bounded, or  $p\geq c(n_{max})^{3r/(r-2)}$ for some positive constant $c$. Thus, by (A.7), (A.8), (A.9) and the control convergence theorem, we can get that
\begin{align*}
\overline{\lim}_{n_{min}\to \infty}\textbf{E}\Big(J_{n,k}^{N \gamma}\Big)&=\overline{\lim}_{n_{min}\to \infty}\textbf{E}\Big(\Big[C_{nk}^{-1} \prod_{i=1}^{k} (Z_{i}Y_i)^{a_{i}}\Big]^{N\gamma}\Big)\leq \Big(1+\gamma \Big)^{\frac{-(k-1)}{2}} \\
\underline{\lim}_{n_{min}\to \infty}\textbf{E}\Big(J_{n,k}^{N \gamma}\Big)&=\underline{\lim}_{n_{min}\to \infty}\textbf{E}\Big(\Big[C_{nk}^{-1} \prod_{i=1}^{k} (Z_{i}Y_i)^{a_{i}}\Big]^{N\gamma}\Big)\geq \Big(1+\gamma \Big)^{\frac{-(k-1)}{2}}.
\end{align*}
for $p\geq c(n_{max})^{3r/(r-2)}$.  This means that 
\begin{align}\tag{A.10}\label{A.10}
\textbf{E}\Big(J_{n,k}^{N \gamma}\Big) \,\,\, \rightarrow \,\,\, \Big(1+\gamma \Big)^{\frac{-(k-1)}{2}} 
\end{align}
as $n_{min} \to \infty$ for $p\geq c(n_{max})^{3r/(r-2)}$. Note that the random variable $ -2log X $ is subject to $ \chi^{2}_{(k-1)} $ distribution if and only if the $ \gamma $-th order moment of the positive random variable $ X $ is equal to $ \Big(1+\gamma \Big)^{\frac{-(k-1)}{2}} $. Thus, let $ X= J_{n,k}^{N} $, by (A.10) we have
\begin{eqnarray*}
 L_{k} = -2 Nlog \, J_{n,k} &=& -2 log\, J_{n,k}^{N} \,\,\, \rightarrow \,\,\, \chi^{2}_{(k-1)},\, as\,\, n_{min}\rightarrow \infty.\,\,\,\,
\end{eqnarray*}
It is clear that for positive number $ \rho_{n} \rightarrow 1 $ we have $ \rho_{n} L_{k}\rightarrow \chi^{2}_{k-1} $ as $ n_{min} \rightarrow \infty $ for fixed $ k\geq 2 $.

\paragraph{•} \textbf{ A.3. Results of  Examples 1 and 2 }

In Examples 1 and 2, we test the null hypothesis $ H_0: \Sigma_{1}  = \Sigma_{2} = \Sigma_{3} = \Sigma$,
where $\Sigma$ is a $ p\times p $ matrix whose $ (a, b) $th element are defined by
$ (-1)^{(a+b)}(0.40)^{\vert a-b\vert} $. We considered the null hypothesis test for two cases of dependent stationary process, the Gaussian and exponential AR(1) models, respectively, with the dimension $ p=350 $.
For both distributions, P-values in table 4 indicate that the $ \rho L_k $ accepts $ H_0 $ when $ H_0 $ is correct. We observe in Table 4 that the $ \rho L_{k} $ test does not reject the null hypothesis $ H_0 $, that is,  each $ P_{value} $ of the $ \rho L_{k} $ test does not fall into the rejection  region $ R_{k} $. The results of  $ \rho L_k $ in Table 4 are completely consistent with the results in Tables 1 and 2.

\begin{table}[H]
\begin{center}
\caption{ The statistic measure and $ P_{value} $ of $ \rho L_{k} $ test for testing a stationary processes AR(1), in Examples 1 and 2.}
    \label{table:table-4}
\begin{tabular}{|c|cc|cc|cc|cc|}\hline
  &\multicolumn{4}{c|}{Example 1.}  &\multicolumn{4}{c|}{Example 2.}  \\\hline
&\multicolumn{2}{c|}{Exponential}  &\multicolumn{2}{c|}{Gaussian}
&\multicolumn{2}{c|}{Exponential}  &\multicolumn{2}{c|}{Gaussian}  \\\hline
    &   $ \chi^{2} $  & $ P_{value} $  & $ \chi^{2} $  & $ P_{value} $
    &   $ \chi^{2} $  & $ P_{value} $  & $ \chi^{2} $  & $ P_{value} $  \\\hline
$ \rho L_{k(\textbf{y}_1)} $ & --- &---- &--- &---- & 1.2006   & 0.7815 & 4.4730  & 0.7094     \\
$ \rho L_{k(\textbf{y}_2)} $ & --- &---- &--- &---- & 1.0615   & 0.5915&  1.3009  & 0.6444    \\
$ \rho L_{k(\textbf{y}_3)}$  & --- &---- &--- &---- & 0.8418   & 0.5267& 2.2446   & 0.8654     \\
$ \rho L_{k(\textbf{y}_4)}$  &--- &---- &--- &---- & 1.2540   & 0.4301 &  2.1401  & 0.8004\\\hline
$ \rho L_{k(\textbf{\textit{1}})}$  &1.0593   & 0.6528 & 1.6076   & 0.7230 & ----   & ---- &  ----  & ----\\
    \hline
\end{tabular}
\end{center}
\end{table}

\end{document}